\documentstyle[12pt]{article}
\catcode`\@=11
\@addtoreset{equation}{section}

\catcode`\@=12
\newtheorem{Theorem}{Theorem}[section]

\newtheorem{Corollary}{Corollary}[section]

\title{Simple-Homotopy Invariance Of Higher Signatures}

\author{Renyi Ma\\
Department of Mathematical Sciences \\
Tsinghua University \\
Beijing, 100084\\
People's Republic of China\\
Email:rma@math.tsinghua.edu.cn}

\date { }

\begin{document}
\textwidth=125mm
\textheight=185mm
\parindent=8mm
\frenchspacing
\maketitle

\begin{abstract}
We prove that the higher signature for any close oriented manifold
is a simple-homotopy invariant.
\end{abstract}

\noindent{\bf 2000 MR Subject Classification} 57N70, 57R95, 57R40

\section{Introduction and results }

Given a closed connected oriented manifold $M^m$ of dimension $m$,
let $[M]$ denote either orientation classes in $H_m(M;Q)$ and in
$H^{m}(M;Z)$, and let ${\cal{L}}_ {M}\in H^{4*}(M;Q)$ be the
Hirzebruch $L$-class of $M$, which is defined as a suitable rational
polynomial in the Pontrjagin classes of $M$ (see\cite{hi}). Denote
the usual Kronecker pairing for $M$, with rational coefficients, by
$$
\left<\,.\,,.\,\right>\colon H^{*}(M;Q)\times
H_{*}(M;Q)\longrightarrow Q\,.
$$

If $M$ is of dimension $m=4k$, then the Hirzebruch Signature Theorem
(see\cite{hi}) says that the rational number
$\left<{\cal{L}}_{M},[M]\right>$ is the signature of the cup product
quadratic form
$$
H^{2k}(M;Z)\otimes H^{2k}(M;Z)\longrightarrow H^{4k}(M;Z)\tilde =Z,
$$
$$
\quad (x,y)\longmapsto (x\cup y)([M])\,.
$$
As a consequence, $\left<{\cal{L}}_{M},[M]\right>$ is an oriented
homotopy invariant of $M$ among close connected oriented manifolds
and of the same dimension $4k$. Let $\alpha\in H^{*}(M;Q)$ be a
prescribed rational cohomology class of $M$. Consider the
\emph{$\alpha$-higher signature}
$$
sign(M,\alpha ) =\big<\alpha\cup{\cal{L}}_{M},[M]\big>\in Q.
$$
Then, the natural conjecture predicts that the rational number
$sign(\alpha ,M)$ is an oriented homotopy invariant of the $M$, in
the sense that $sign(\alpha ,M)=sign(\beta  ,N)$ whenever $N^{m}$ is
a second close connected oriented manifold and there exists a
homotopy equivalence $h:M\to N$ such that $h^*\beta=\alpha$.
\begin{Theorem}
the higher signature $sign(M,\alpha )$ coming from $(M,\alpha )$ are
oriented simple-homotopy invariants for closed connected oriented
manifolds of arbitrary dimension.
\end{Theorem}

\begin{Corollary}
the famous Novikov Conjecture holds for the close manifold $M$ with
trivial Whitehead group (see\cite{gro-mm}).
\end{Corollary}

\section{Proof of Theorem1.1}

\vskip 3pt

 {\bf Proof of Theorem1.1:} Let $M$ be a close oriented connected manifold.
Let $M'$ be a second close connected oriented manifold which is
simple-homotopy equivalent to $M$ by the maps $h:M\to M'$ and
$h':M'\to M$, here $h'\circ h \cong _sId_M:M\to M $ and $h\circ h
'\cong _sId_{M'}:M'\to M'$. So, by Mazur's theorem, i.e., the
non-stable neighborhood theorem(see\cite{ma}), we have a
diffeomorphism
$$
f:M\times R^{m+1}\to M'\times R^{m+1}.
$$
So,
$$
{{\cal {L}}}_M={\cal {L}}_{(M\times R^{m+1})}={\cal {L}}_{(M'\times
R^{m+1})}={\cal {L}}_{M'}.
$$
This yields Theorem1.1.

\end{document}